\documentclass[10pt]{article}

\usepackage[dvips]{epsfig}
\usepackage{latexsym}
\usepackage{amssymb}
\usepackage{algorithm}
\usepackage{algorithmic}
\usepackage{amsfonts}
\usepackage{amsmath}
\usepackage{indentfirst}

\title{A closed formula for the number of convex permutominoes}

\newtheorem{proposition}{Proposition}

\newcommand{\qed}{\hfill\square\medskip}

\author{
 Filippo Disanto
 \thanks{Universit\`a di Siena, Dipartimento di Scienze Matematiche e Informatiche, Pian dei Mantellini
 44, 53100 Siena, Italy ({\tt rinaldi@unisi.it}).}
 \and
 Andrea Frosini
 \thanks{Universit\`a di Firenze, Dipartimento di Sistemi e Informatica, viale Morgagni 65, 50134 Firenze, Italy
 ({\tt [frosini,  pinzani]@dsi.unifi.it}).}
 \and
 Renzo Pinzani $^\dag$
 \and
 Simone Rinaldi
 $^*$
}

\date{} 

\begin{document}\date{}

\maketitle

\begin{abstract}
In this paper we determine a closed formula for the number of convex permutominoes of size $n$. We reach this
goal by providing a recursive generation of all convex permutominoes of size $n+1$ from the objects of size $n$,
according to the ECO method, and then translating this construction into a system of functional equations
satisfied by the generating function of convex permutominoes. As a consequence we easily obtain also the
enumeration of some classes of convex polyominoes, including stack and directed convex permutominoes.
\end{abstract}

\section{Basic definitions and contents of the paper}

\begin{figure}[htb]
\begin{center}
\epsfig{file=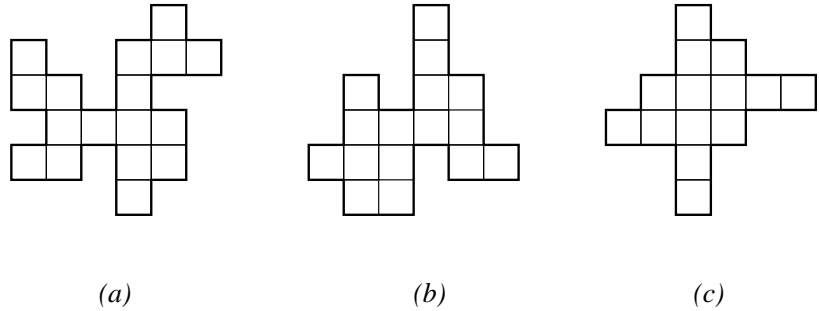}
\end{center}
\caption{(a) a polyomino; (b) a column convex polyomino; (c) a convex polyomino.
\label{polyom}}
\end{figure}

A {\em polyomino} is a finite union of elementary cells of  the lattice $Z \times Z$, whose interior is connected
(see Figure~\ref{polyom}~$(a)$). Polyominoes are defined up to translations. A polyomino is said to be {\em
column convex} (resp. {\em row convex}) if every its column (resp. row) is connected (see
Figure~\ref{polyom}~$(b)$). A polyomino is said to be {\em convex}, if it is both row and column convex (see
Figure~\ref{polyom}~$(c)$).

\medskip

The number $f_{n}$ of convex polyominoes with semi-perimeter $n+2$ was determined by Delest and Viennot,
in \cite{DV}:

\begin{equation}\label{eq1}
c_{n+2}=(2n+11)4^n \, - \, 4(2n+1) \, {2n \choose n}, \quad n \geq 0; \quad c_{0}=1, \quad c_{1}=2,
\end{equation}

\noindent sequence $A005436$ in \cite{sloane}, the first few terms being:
$$1,2,7,28,120,528,2344,10416, \ldots .$$

In the last two decades convex polyominoes, and several combinatorial objects obtained as a generalizations of
this class, have been studied by various points of view. For the main results concerning the enumeration and
other combinatorial properties of convex polyominoes we refer to \cite{mbm2, mbm, gutman, chang}.

\subsection{Permutominoes}

Let $P$ be a polyomino without holes, having $n$ rows and columns, $n\geq 1$;
we assume without loss of generality that the south-west
corner of its minimal bounding rectangle is placed in $(1,1)$.
Let ${\cal A}= \{ A_1, \ldots ,A_{2(r+1)} \}$ be the set of its
vertices ordered in a clockwise sense starting from the leftmost vertex
having minimal ordinate.

We say that $P$ is a {\em permutomino} if the sets ${\cal P}_1 = \{ A_1, A_3, \ldots ,A_{2r+1} \}$ and ${\cal
P}_2 = \{ A_2, A_4, \ldots ,A_{2r+2} \}$ represent two permutation matrices of $[n+1]=\{ 1, 2, \ldots ,n+1\}$.
Obviously, if $P$ is a permutomino, then $r=n$, and $n$ is called the {\em size} of the permutomino.

\begin{figure}[htb]
\begin{center}
\epsfig{file=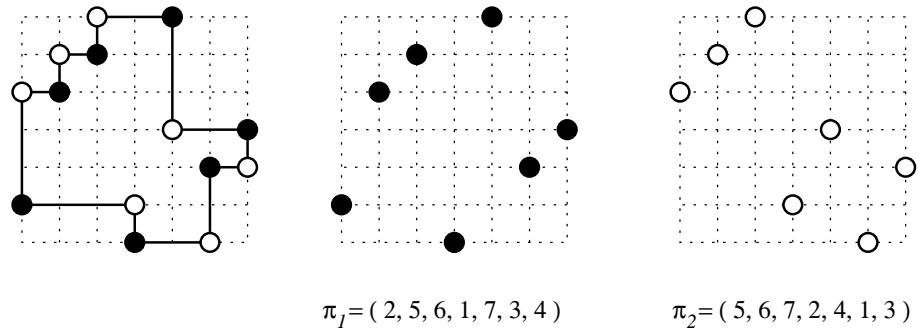} \caption{A permutomino and the two associated permutations. \label{permu1}}
\end{center}
\end{figure}

\begin{figure}[htb]
\begin{center}
\epsfig{file=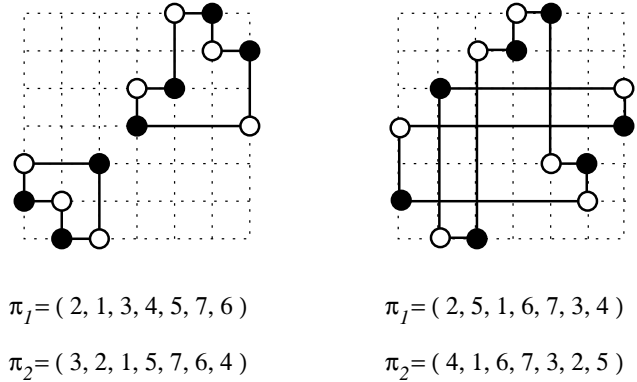} \caption{The two main cases when a couple of permutations $\pi _1$ and $\pi _2$ of $[n]$,
with $\pi_1 (i) \neq \pi_2 (i)$, $1 \leq i \leq n+1$, does not define a permutomino: (a) two disconnected sets of
cells; (b) the boundary intersects itself. \label{permuz}}
\end{center}
\end{figure}

The two permutations associated with ${\cal P}_1$ and ${\cal P}_2$ are indicated by $\pi _1$ and $\pi _2$,
respectively (see Figure \ref{permu1}). While it is clear that any permutomino of size $n$ uniquely individuates
two point-by-point distinct permutations $\pi _1$ and $\pi _2$ of $[n+1]$, not all the couples of permutations
$\pi _1$ and $\pi _2$ of $n$ such that $\pi_1 (i) \neq \pi_2 (i)$, $1 \leq i \leq n+1$ define a permutomino, as
it was partially investigated in \cite{fanti}  (see Figure \ref{permuz}).

\medskip

Permutominoes were introduced by F.~Incitti in \cite{incitti} while studying the problem of determining the
$\widetilde{R}$-polynomials (related with the Kazhdan-Lusztig R-polynomials) associated with a pair $(x,y)$ of
permutations. Concerning the class of polyominoes without holes, our definition (though different) turns out to
be equivalent to Incitti's one, which is more general but uses some algebraic notions not necessary in this
paper.

\medskip

In this paper we deal with the enumeration of convex polyominoes which are also permutominoes, the so called {\em
convex permutominoes}. From the definition we have that in any convex permutomino $P$, for each abscissa
(ordinate) there is exactly one vertical (horizontal) side in the boundary of $P$ with that coordinate. It is
simple to observe that the previous property is also a sufficient condition for a convex polyomino to be a
permutomino.

In \cite{fanti}, using bijective techniques, it was proved that the number of {\em parallelogram permutominoes}
of size $n$ is equal to the $n$th {\em Catalan number},

$$ \frac{1}{n+1} \, {{2n} \choose {n}}, $$

\noindent and moreover, that the number of {\em directed-convex permutominoes} of size $n$ is equal to half the
$n$th {\em binomial coefficient},

$$ \frac{1}{2} \, {{2n} \choose {n}}. $$

\medskip

The first attempt to count convex permutominoes was made in \cite{milanesi}, where the authors considered the
couples of permutations which define convex permutominoes, and then obtained an expression for the number $f_n$
of convex permutominoes of size $n$:

$$ f_{n+1} = \sum_{s=0}^{n-2} \sum_{t=0}^s \sum_{x=0}^t \, {{n-2} \choose t} \, {{n-2} \choose t} \, {{n-2} \choose {x+s-t}}
\, - \, (n-1) \, {{2(n-2)} \choose {n-2}} \, + \, 4^{n-2} ,$$

\noindent but then were not able to derive the generating function, nor the closed form for the number $f_n$.

In this paper we deal with the same enumeration problem using a different and more immediate approach: we
determine a direct recursive construction for the convex permutominoes of a given size, based on the application
of the ECO method, which easily leads to the generating function, and finally prove that the number of convex
permutominoes of size $n$ is:

\begin{equation}
    f_n \, = \, 2 \, (n+3) \, 4^{n-2} \, - \, \frac{n}{2} \, {{2n} \choose {n}} \qquad n\geq 1.
\end{equation}

\subsection{ECO method}

In this section we will recall some basics about the ECO method, where ECO stands for Enumeration of
Combinatorial Objects. Such a method, introduced by Pinzani and his collaborators in \cite{Eco1}, is a
constructive method to produce all the objects of a given class, according to the growth of a certain parameter
(the \emph{size}) of the objects. Basically, the idea is to perform ``local expansions'' on each object of size
$n$, thus constructing a set of objects of the successive size (see \cite{Eco1} for more details).

The application of the ECO method often leads to an easy solution for problems that are commonly believed
``hard'' to solve. For example, in \cite{DDFR} the authors give an ECO construction for the classes of convex
polyominoes and column-convex polyominoes according to the semi-perimeter. A simple algebraic computation leads
then to the determination of generating functions for the two classes.

In \cite{exh} it is also shown that an ECO construction easily leads to an efficient algorithm for the exhaustive
generation of the examined class. Moreover, an ECO construction can often produce interesting combinatorial
information about the class of objects studied, as shown in \cite{Eco1} using analytic methods, or in \cite{BFR},
using bijective techniques. In \cite{GFGT}, Banderier et al. reintroduced the  {\em kernel method} in order to
determine the generating function of various types of ECO systems.

\medskip

Going deeper into formalism, let $p$ be a parameter $p:{\cal O}\to \mathbb{N}^+$, such that $\left | {\cal O}_n
\right | = \left | \, \{ O\in {\cal O}:p(O)=n \} \, \right |$ is finite. An operator $\vartheta$ on the class
${\cal O}$ is a function from ${\cal O}_n$ to $2^{{\cal O}_{n+1}}$, where $2^{{\cal O}_{n+1}}$ is the power set
of ${{\cal O}_{n+1}}$.

\begin{proposition}\label{eco}
{\em Let $\vartheta$ be an operator on $\cal O$. If $\vartheta$ satisfies the following conditions:

\begin{description}
\item[1.]  for each $O' \in {\cal O}_{n+1}$, there exists $O\in
{\cal O}_n$ such that $O'\in \vartheta (O)$, \item[2.] for each
$O,O'\in {\cal O}_n$ such that $O\neq O'$, then $\vartheta (O)
\cap \vartheta(O') = \emptyset $,\end{description}

\noindent then the family of sets ${\cal F}_{n+1}=\{ \vartheta (O): O\in {\cal O}_n \}$ is a partition of ${\cal
O}_{n+1}$.}
\end{proposition}

This method was successfully applied to the enumeration of various
classes of walks, permutations, and polyominoes. We refer to
\cite{Eco1}, and \cite{rule} for further details and results.

The recursive construction determined by $\vartheta$ can be
suitably described through a {\em generating tree}, i.e. a rooted
tree whose vertices are objects of $\cal O$. The objects having
the same value of the parameter $p$ lie at the same level, and the
sons of an object are the objects it produces through $\vartheta$.

If the construction determined by the ECO operator $\vartheta$ is
regular enough it is then possible to describe it by means of a
{\em succession rule} of the form:

$$
\left\{
\begin{array}{l}
(b) \\
 \\
(h) \rightsquigarrow (c_1)(c_2)\ldots (c_{q(h)}),
\end{array}
\right.
$$

\noindent where $b, h, c_i \in \Bbb N$, and $q: \Bbb N^+ \to \Bbb N^+$, meaning that object at the root of the
generating tree has $b$ sons, and the $q(h)$ objects $O'_1, \ldots ,O'_{q(h)}$, produced by an object $O$ are
such that $\left | \vartheta (O'_i) \right | = c_i$, $1 \leq i \leq q(h)$. A succession rule defines a sequence
$\{f_n\}_{n\geq 1}$ of positive integers, where $f_n$ is the number of nodes at level $n$ of the generating tree,
assuming that the root is at level $1$.

\section{Generation of convex permutominoes}

Let ${\cal C}_n$ be the set of convex permutominoes of size $n$. In order to define the ECO construction for
convex permutominoes, we need to point out a simple property of their boundary, related to {\em reentrant} and
{\em salient points}. So let us briefly recall the definition of these objects.

\begin{figure}[htb]
\begin{center}
\epsfig{file=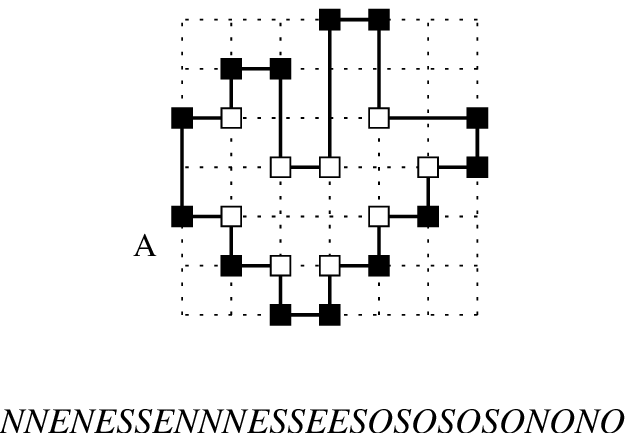} \caption{The coding of the boundary of a polyomino, starting from $A$ and moving in a
clockwise sense; its salient (resp. reentrant) points have been evidenced by a black (resp. white) square.
\label{reentrant}}
\end{center}
\end{figure}

Let $P$ be a polyomino; starting from the leftmost point having minimal ordinate, and moving in a clockwise
sense, the boundary of $P$ can be encoded as a word in a four letter alphabet, $\{ N, E, S, W \}$, where $N$
(resp. $E$, $S$, $W$) represents a {\em north} (resp.{\em east}, {\em south}, {\em west}) unit step. Any
occurrence of a sequence $NE$, $ES$, $SO$, or $ON$ in the word encoding $P$ defines a {\em salient point} of $P$,
while any occurrence of a sequence $EN$, $SE$, $OS$, or $NO$ defines a {\em reentrant point} of $P$ (see for
instance, Figure \ref{reentrant}).

In \cite{daurat} and successively in \cite{brlek}, in a more general context, it was proved that in any polyomino
the difference between the number of salient and reentrant points is equal to $4$.

Let us turn to consider the class of convex permutominoes. In a convex permutomino of size $n$ the length of the
word coding the boundary is $4n$, and we have $n+3$ salient points and $n-1$ reentrant points; moreover we
observe that a reentrant point cannot lie on the minimal bounding rectangle. This leads to the following
remarkable property:

\begin{proposition}\label{spigoli}
The set of reentrant points of a convex permutomino of size $n$ defines a permutation matrix of $[n-1]$, $n\geq
2$.
\end{proposition}

For simplicity of notation, and to clarify the definition of the upcoming ECO construction, we agree to group the
reentrant points of a convex permutomino in four classes; in practice we choose to represent the reentrant point
determined by a sequence $EN$ (resp. $SE$, $OS$, $NO$) with the symbol $\alpha$ (resp. $\beta$, $\gamma$,
$\delta$). Using this notation we can state that each convex permutomino of size $n\geq 2$ can be uniquely
represented by the permutation matrix defined by its reentrant points, which has dimension $n-1$, and uses the
symbols $\alpha , \, \beta , \, \gamma , \, \delta $.

\bigskip

\begin{figure}[htb]
\begin{center}
\epsfig{file=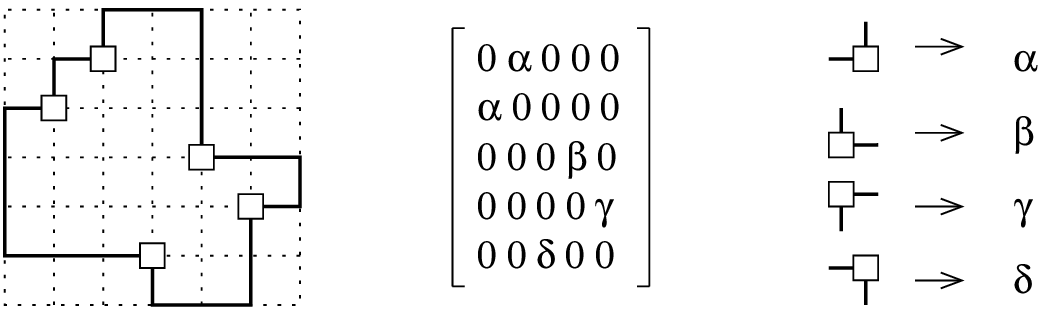} \caption{The reentrant points of a convex permutomino uniquely define a permutation
matrix in the symbols $\alpha$, $\beta$, $\gamma$ and $\delta$. \label{permu2}}
\end{center}
\end{figure}

\subsection{The ECO operator}

Let $P\in {\cal C}_n$; the number of cells in the
rightmost column of $P$ is called the {\em degree} of $P$.
For any $n\geq 2$ we partition the class ${\cal C}_n$
into three distinct classes. In order to define these classes, let us consider the following
conditions on a convex permutomino:

\begin{description}
\item[U1]: the uppermost cell of the rightmost column of the
polyomino has the maximal ordinate among all the cells of the
polyomino; \item[U2]: the lowest cell of the rightmost column of
the polyomino has the minimal ordinate among all the cells of the
polyomino.
\end{description}

\begin{figure}[htb]
\begin{center}
\epsfig{file=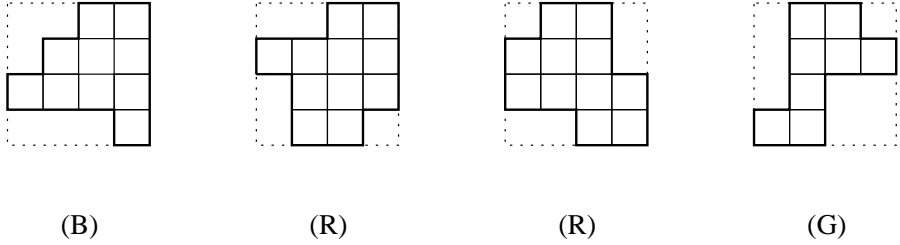,width=4.7in,clip=} \caption{Convex permutominoes in classes $B$, $R$, and $G$.
\label{4classes}}
\end{center}
\end{figure}

\noindent We say that a convex permutomino $P$ belongs to class:

\begin{description}
    \item{-} $B$, if it satisfies both conditions {\bf U1} and {\bf U2} (i.e. $P$ has degree $n$,
    see Figure \ref{4classes}, (B));
    \item{-} $R$, if it satisfies only one among conditions {\bf U1}, {\bf U2} (see Figure \ref{4classes}, (R));
    \item{-} $G$, if it satisfies none of conditions {\bf U1}, {\bf U2} (see Figure \ref{4classes}, (G)).
\end{description}

For simplicity sake, each permutomino in class $B$ (resp. $R$, $G$) and degree $k$ is represented by the label
$(k)_b$ (resp. $(k)_r$, $(k)_g$). For instance, the four permutominoes depicted in Figure \ref{4classes} have
labels $(4)_b$, $(3)_r$, $(2)_r$, $(1)_g$, respectively. We assume that the single cell permutomino belongs to
class $B$, then it has the label $(1)_b$.

\bigskip

Our aim is now to use the property stated in Proposition \ref{spigoli} to define an ECO operator $\vartheta :
{\cal C}_n \to 2^{{\cal C}_{n+1}}$ which defines a recursive construction of all the objects of size $n+1$ in a
unique way from the objects of size $n$. The operator $\vartheta$ acts on a convex permutomino performing some
local expansions on the cells of its rightmost column. In order to define these operations let us consider a
generic permutomino $P$ of size $n$, let us indicate by $c_1, \ldots ,c_n$ (resp. $r_1, \ldots ,r_n$) the columns
(resp. rows) of $P$ numbered from left to right (resp. bottom to top), and by $\ell (c_i)$ (resp. $\ell (r_i)$)
the number of cells in the $i$th column (resp. $i$th row), with $1 \leq i \leq n$. The four operations of
$\vartheta$ will be denoted by $\alpha$, $\beta$, $\gamma$, and $\delta$, and below we give a detailed
description of each of them:

\begin{figure}[htb]
\begin{center}
\epsfig{file=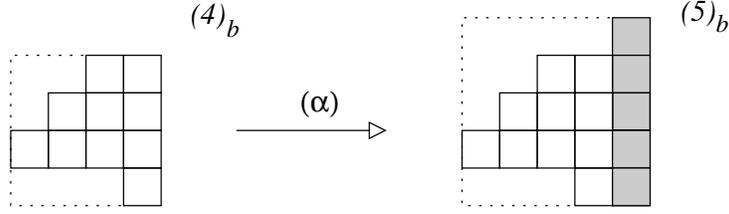} \caption{Operation $(\alpha )$ performed on a permutomino of class $B$. The added
column has been highlighted. \label{operationA}}
\end{center}
\end{figure}

\begin{description}
  \item{$( \alpha )$} if $P$ satisfies condition {\bf U1}, then $( \alpha )$ adds a new column made of $c_n+1$ cells
  on the right of $c_n$, according to Figure \ref{operationA}.

It is clear that the obtained polyomino is a convex permutomino of size $n+1$, still satisfying condition {\bf
U1}; the rightmost reentrant point in such new permutomino is of type $\alpha$ (this is the reason why we have
called the reentrant points with the same name of the operations on permutominoes).

\item{$(\beta )$} it can be performed on each cell of $c_n$; so let $d_i$ be the $i$th cell of $c_n$, from bottom
to top, with $1 \leq i \leq \ell (c_n)$. Operation $( \beta )$ adds a new row above the row containing $d_i$ (of
the same length), and a new column on the right of $c_n$, made of $i$ cells, as illustrated in Figure
\ref{operationB}.

Observe that, since the new added row  is long as the row below it, we ensure that the obtained polyomino has a
unique horizontal side at level $i$, while adding the new column from bottom to level $i$ we ensure that the
obtained polyomino has a unique vertical side at abscissa $n-1$, hence the basic property of permutominoes is
preserved.

\begin{figure}[htb]
\begin{center}
\epsfig{file=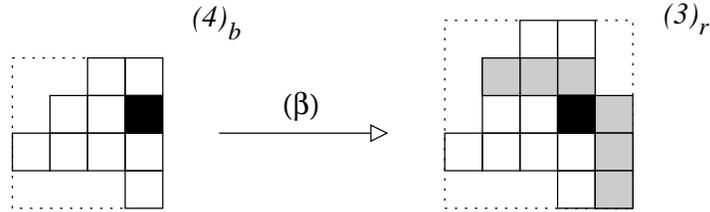} \caption{Operation $(\beta )$ performed on a cell $d_i$ of the rightmost column of a
polyomino in class $B$. The cell $d_i$ is filled in black, the added row and column have been highlighted.
\label{operationB}}
\end{center}
\end{figure}

Then it is clear that, for any $i$, the obtained polyomino is a convex permutomino of size $n+1$, and its
rightmost reentrant point is of type $\beta$.

\item{$( \gamma )$} it can be performed on each cell of $c_n$; so let $d_i$ be the $i$th cell of $c_n$, from
bottom to top, with $1 \leq i \leq \ell (c_i)$. Operation $( \gamma )$ adds a new row below the row containing
$d_i$ (of the same length), and a new column on the right of $c_n$, made of $n-i+1$ cells, as illustrated in
Figure \ref{operationC}.

\begin{figure}[htb]
\begin{center}
\epsfig{file=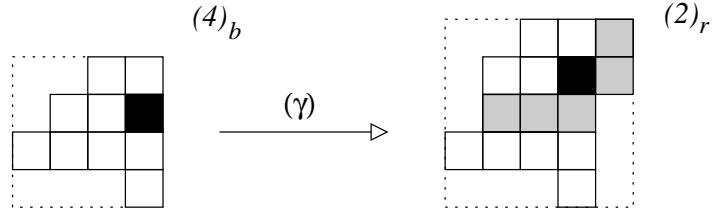} \caption{Operation $( \gamma )$ performed on a cell $d_i$ of the rightmost column of
a polyomino in class $B$. The cell $d_i$ is filled in black, the added row and column have been highlighted.
\label{operationC}}
\end{center}
\end{figure}

It is clear that, for any $i$, the obtained polyomino is a convex permutomino of size $n+1$, and its rightmost
reentrant point is of type $\gamma$.

\item{$( \delta )$} if $P$ satisfies condition {\bf U2}, then $(\delta )$ adds a new column made of $c_n+1$ cells
  on the right of $c_n$, according to Figure \ref{operationD}.

It is clear that the obtained polyomino is a convex permutomino of size $n+1$, still satisfying condition {\bf
U2}; the rightmost reentrant point in such new permutomino is of type $\delta$.

\begin{figure}[htb]
\begin{center}
\epsfig{file=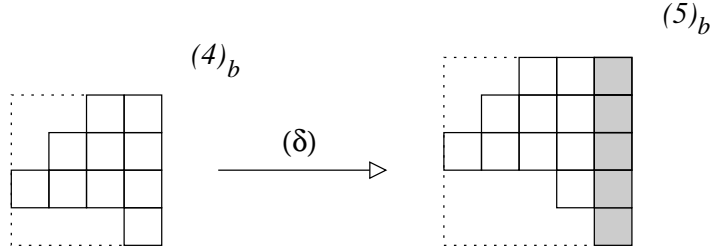} \caption{Operation $(\delta )$ performed a polyomino in class $B$.
\label{operationD}}
\end{center}
\end{figure}

\end{description}

As we already mentioned, the operations performed by $\vartheta$ on a convex permutomino $P$
depend on the family to which $P$ belongs. So let us consider the different cases:

\begin{enumerate}
    \item $P$ {\em belongs to class $B$}. The operator $\vartheta$ performs on $P$ operations
    $( \alpha )$, $( \delta )$ and one application of $(\beta )$ and $(\gamma )$ for any cell in $c_n$.
    So, let $k$ be the degree of $P$, the application of $\vartheta$ to $P$ produces
    $2k+2$ different convex permutominoes of size $n+1$ (see Figure \ref{esempioB}).

\begin{figure}[htb]
\begin{center}
\epsfig{file=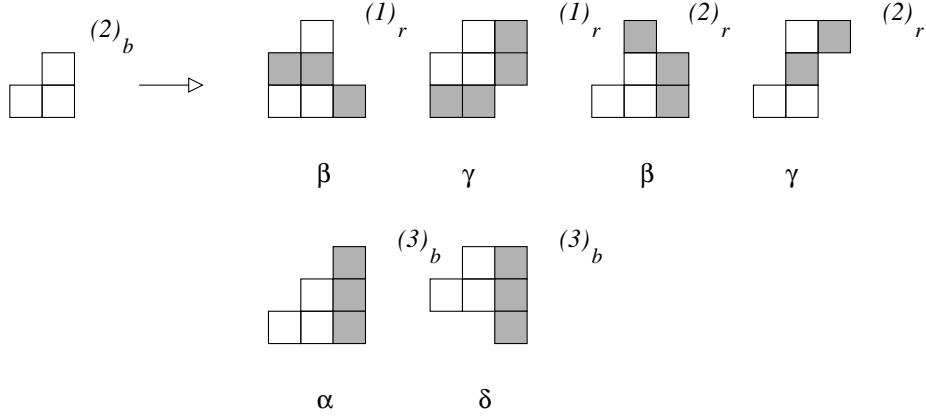,width=4.8in,clip=} \caption{The operator $\vartheta$ applied to a permutomino of class $B$;
the added rows and columns are highlighted,
and the applied operation is mentioned below. \label{esempioB}}
\end{center}
\end{figure}

    More formally, applying $\vartheta$
    to a convex permutomino of label $(k)_b$, we have $2(k+1)$ different permutominoes,
    two with label $(1)_r, (2)_r, \ldots (k)_r$, and two with label $(k+1)_b$.
    This can be formalized by the production:

    $$ (k)_b \rightsquigarrow (1)_r\, (1)_r \, (2)_r \, (2)_r \, \ldots \, (k)_r \, (k)_r \, (k+1)_b \,  (k+1)_b .$$

    \item {\em $P$ belongs to class $R$}. There are two possibilities:

    \begin{description}
        \item{i.} $P$ satisfies {\bf U1} (and not {\bf U2}).
        The operator $\vartheta$ performs on $P$ operation
        $( \alpha )$, and one application of operations $(\beta )$ and $(\gamma )$ for any cell in $c_n$.

        \item{ii.} $P$ satisfies {\bf U2} (and not {\bf U1}). The operator $\vartheta$ performs the following operations:
        The operator $\vartheta$ performs on $P$ operations
        $( \delta )$, and one application of operations $(\beta )$ and $(\gamma )$ for any cell in $c_n$ (see Figure~\ref{esempio}).

\begin{figure}[htb]
\begin{center}
\epsfig{file=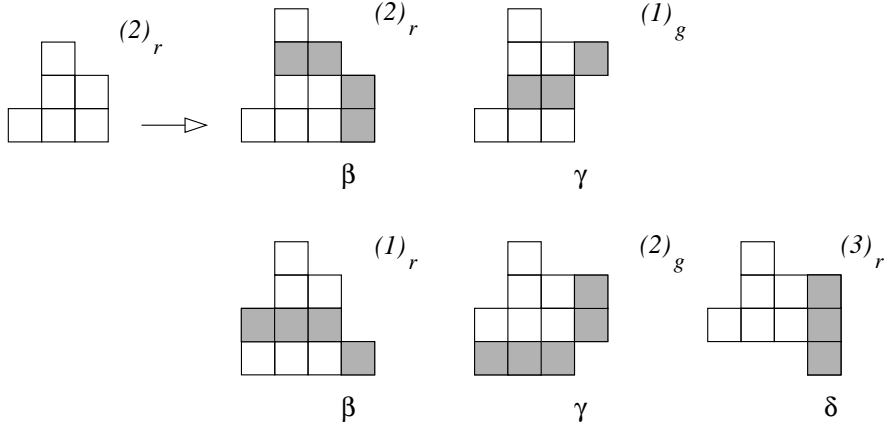,width=4.6in,clip=} \caption{The operator $\vartheta$ applied to a permutomino of class
$R$, satisfying {\bf U2} (and not {\bf U1}); the added rows and columns are highlighted, and the applied
operation is mentioned below. \label{esempio}}
\end{center}
\end{figure}

    \end{description}

        In both cases, being $k$ be the degree of $P$, the application of $\vartheta$ to $P$ produces
        $2k+1$ different convex permutominoes of size $n+1$.
        More formally, applying $\vartheta$
        to a convex permutomino of label $(k)_r$, we have $2k+1$ different permutominoes,
        with labels $(1)_r, (2)_r, \ldots (k)_r, (k+1)_r$, and $(1)_g, (2)_g, \ldots (k)_g$.
        This can be formalized by the production:

    $$ (k)_r \rightsquigarrow (1)_r\, (1)_g \, (2)_r \, (2)_g \, \ldots \, (k)_r \, (k)_g \, (k+1)_r.$$

\item $P$ {\em belongs to class} $G$. The operator $\vartheta$ performs on $P$ an application of
    operations $(\beta )$ and $(\gamma )$ for any cell in $c_n$.
    So, let $k$ be the degree of $P$, the application of $\vartheta$ to $P$ produces
    $2k$ different convex permutominoes of size $n+1$.
    More formally, applying $\vartheta$
    to a convex permutomino of label $(k)_g$, we have $2k$ different permutominoes,
    two with labels $(1)_g, (2)_g, \ldots (k)_g$.
    This can be formalized by the production:

    $$ (k)_g \rightsquigarrow (1)_g\, (1)_g \, (2)_g \, (2)_g \, \ldots \, (k)_g \, (k)_g .$$

\end{enumerate}

\bigskip

\begin{proposition}
The operator $\vartheta$ satisfies conditions 1. and 2. of Proposition \ref{eco}.
\end{proposition}

\noindent {\bf Proof.} We have to prove that any convex permutomino of size $n\geq 2$ is uniquely obtained
through the application of the operator $\vartheta$ to a convex permutomino of size $n-1$. So let $P\in {\cal
C}_{n}$, and, as usual, let us indicate by $c_1, \ldots ,c_n$ (resp. $r_1, \ldots ,r_n$) the columns (resp. rows)
of $P$ numbered from left to right (resp. bottom to top), and by $\ell (c_i)$ (resp. $\ell (r_i)$) the number of
cells in the $i$th column (resp. $i$th row), with $1 \leq i \leq n$. We look at the rightmost reentrant point of
$P$, which is unique due to Proposition \ref{spigoli}, and we have the following four possibilities:

\begin{enumerate}
    \item the rightmost reentrant point of $P$ is of type $\alpha$, i.e. $\ell (c_n)=n$; due to the permutomino
    definition, it is clear that $\ell (r_n)=1$, then $P$ has been produced through the application of operation
        $( \alpha )$ to the permutomino $P'\in {\cal C}_{n-1}$, obtained removing column $c_n$ from $P$
        (see Figure~\ref{operationA});

    \item the rightmost reentrant point of $P$ is of type $\beta$, and then necessarily $1 \leq \ell (c_n)<n$;
    let $P'$ be the permutomino of ${\cal C}_{n-1}$ obtained by removing the column $c_n$ and the row $r_{\ell
    (c_n)+1}$ from $P$. Is is then clear that $P$ is produced through the application of operation $(\beta )$ to
    the $\ell (c_n)$-th cell (from bottom to top) of $P'$ (see Figure~\ref{operationB});

    \item the rightmost reentrant point of $P$ is of type $\gamma$, and then necessarily $1 \leq \ell (c_n)<n$;
    let $P'$ be the permutomino of ${\cal C}_{n-1}$ obtained by removing the column $c_n$ and the row $r_{n - \ell
    (c_n)}$ from $P$. Is is then clear that $P$ is produced through the application of operation $(\gamma )$ to
    the $(n - \ell (c_n)+1)$-th cell (from bottom to top) of $P'$ (see Figure~\ref{operationC});

  \item the rightmost reentrant point of $P$ is of type $\delta$, also with $\ell (c_n)=n$; due to the permutomino
    definition, it is clear that $\ell (r_1)=1$, then $P$ has been produced through the application of operation
    $(\delta )$ to the permutomino $P'\in {\cal C}_{n-1}$, obtained removing column $c_n$ from $P$
    (see Figure~\ref{operationD}). $\qed$

\end{enumerate}

\medskip

The growth of convex permutominoes defined by the ECO operator $\vartheta$ can be suitably represented in terms
of the succession rule $\Omega$:

$$
\Omega: \left\{\begin{array}{l}
(1)_b \\
\\
(k)_b \rightsquigarrow (1)_r \, (1)_r \, \ldots \, (k)_r \, (k)_r
\, (k+1)_b \, (k+1)_b \\
\\
(k)_r \rightsquigarrow (1)_r \, (1)_g \, \ldots \, (k)_r \, (k)_g
\, (k+1)_r \\
\\
(k)_g \rightsquigarrow (1)_g \, (1)_g \, \ldots \, (k)_g \, (k)_g
\, .
\end{array}
\right.
$$

\begin{figure}
\begin{center}
\epsfig{file=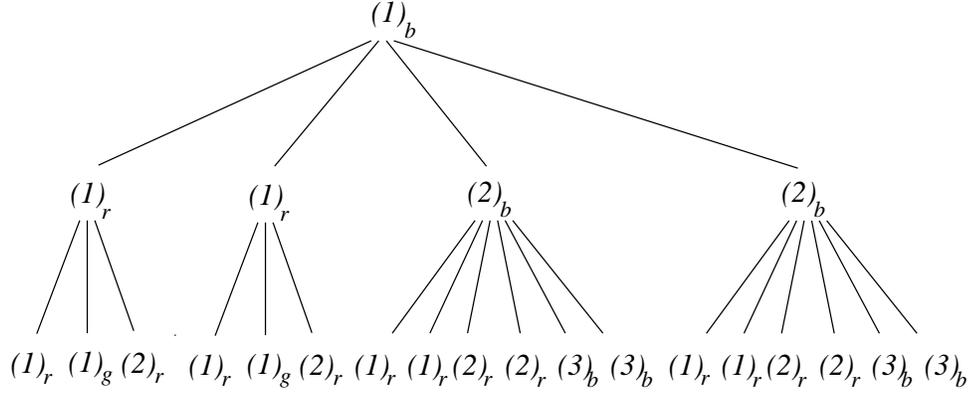} \caption{The first three levels of the
generating tree of the rule $\Omega$. \label{bij2}}
\end{center}
\end{figure}

\bigskip

The root of the tree is $(1)_b$, which is the label of the one cell polyomino.

\bigskip

\section{Enumeration of convex permutominoes}

In this section we will determine the generating function of
convex permutominoes according to various parameters, using the
simple remark that the number $f_n$ of convex permutominoes of
size $n$ is given by the number of objects at level $n$ of the
generating tree of $\Omega$, $n \geq 1$, assuming without loss of generality that
the root of the tree is at level $1$.

\medskip

\noindent Throughout this section, we will use the following notation:

\begin{description}
\item{-} $\mathrm{F}$ is the set of labels of the generating tree
of $\Omega$;
\item{-} $\mathrm{B}$ (resp. $\mathrm{R}$, $\mathrm{G}$) is the set of labels $(k)_b$ (resp. $(k)_r$,
$(k)_g$), $k \geq 1$, in the generating tree of $\Omega$.
\end{description}

\noindent Moreover, for any convex permutomino $P$, let $l(P)$ (briefly,
$l$) be the label of $P$, and $n(P)$ (briefly,
$n$) be the size of $P$. Our aim is to determine the
generating function:

$$ F(s,t)=\sum _{P\in {\mathrm{F}}} s^{l(P)}t^{n(P)}= st+(2s+2s^2)t^2+(8s+6s^2+4s^3)t^3 + \ldots, $$

\noindent since in particular $F(1,t)$ is the generating function of convex permutominoes according to the size.
To do this we need to consider the following auxiliary generating functions:

\begin{description}
    \item[-] $B(s,t)=\sum _{P\in {\mathrm{B}}} s^{l(P)}t^{n(P)}$, i.e. the generating function of $\mathrm{B}$,
    \item[-] $R(s,t)=\sum _{P\in {\mathrm{R}}} s^{l(P)}t^{n(P)}$, i.e. the generating function of $\mathrm{R}$,
    \item[-] $G(s,t)=\sum _{P\in {\mathrm{G}}} s^{l(P)}t^{n(P)}$, i.e. the generating function of $\mathrm{G}$.
\end{description}

\noindent Clearly, $F(s,t)=B(s,t)+R(s,t)+G(s,t)$. From the
productions of $\Omega$ we obtain the following relations concerning $B(s,t)$:

$$ B(s,t)= st + \sum_{P \in {\mathrm{B}}} 2 \, s^{l+1} \, t^{n+1}, $$

\noindent hence we have that

\begin{equation}\label{bst}
B(s,t)=  \frac{st}{1-2st} \qquad B(1,t)=  \frac{t}{1-2t}
\end{equation}

\noindent Then we are able to write down the equation for $R(s,t)$:

\begin{eqnarray*}
  R(s,t) &=&  \sum_{P \in {\mathrm{B}}} 2 \, \left ( s \, t^{n+1} + \ldots + s^l t^{n+1} \right ) \, +
  \sum_{P \in {\mathrm{R}}} \left ( s \, t^{n+1} + \ldots + s^{l+1} t^{n+1} \right ) \\
    &=& 2st \sum_{P \in {\mathrm{B}}} \, \frac{1-s^l}{1-s} \, t^n \, + \,
    st \sum_{P \in {\mathrm{R}}} \, \frac{1-s^{l+1}}{1-s} \, t^n \\
    &=& \frac{2st}{1-s} \, \left ( B(1,t) - B(s,t) \right ) \, + \, \frac{st}{1-s} \, \left ( R(1,t) - s \, R(s,t) \right ).
\end{eqnarray*}

\noindent Let $\overline{B}(s,t)=B(1,t)-B(s,t)$; the previous equation can be re-written as:

\begin{equation}\label{rst}
R(s,t) \, \left ( 1+ \frac{s^2t}{1-s} \right) \, = \, \frac{2st}{1-s} \, \overline{B}(s,t) \, + \,
\frac{st}{1-s} \, R(1,t).
\end{equation}

\noindent Equation (\ref{rst}) has two unknowns: $R(s,t)$, and $R(1,t)$. Applying the {\em kernel method} (as explained in
detail in \cite{GFGT}) we look for the value of $s$ for which the factor on the left multiplying $R(s,t)$
is equal to zero, i.e. the solution of the kernel

$$ 1-s+ts^2=0 .$$

\noindent Of the two solutions we observe that only

$$ s_0=\frac{1-\sqrt{1-4t}}{2t} $$

\noindent is a formal power series with positive coefficients. Substituting $s=s_0$ in $(\ref{rst})$ we have:

$$ \frac{2s_0t}{1-s_0} \, \overline{B}(s_0,t) \, + \,
\frac{s_0t}{1-s_0} \, R(1,t) \, = \, 0, $$

\noindent which leads to:

\begin{equation}\label{r1t}
R(1,t)=\frac {1}{\sqrt {1-4\,t}} - { \frac{1}{1-2\,t}}.
\end{equation}

\noindent Replacing the value of $R(1,t)$ into (\ref{rst}) we obtain

\begin{equation}\label{r1t}
R(s,t)=
{\frac {ts \left( 2\,ts-1+\sqrt {1-4\,t} \right) }{\sqrt {1-4\,t}
 \left( 2\,{t}^{2}{s}^{3}+2\,ts-3\,t{s}^{2}-1+s \right) }}.
\end{equation}

\noindent Finally, from the productions of $\Omega$ we derive the following equation for $N(s,t)$:

\begin{eqnarray*}
  N(s,t) &=&  \sum_{P \in {\mathrm{R}}} \, \left ( s \, t^{n+1} + \ldots + s^l t^{n+1} \right ) \, +
  \sum_{P \in {\mathrm{N}}} \, 2 \, \left ( s \, t^{n+1} + \ldots + s^{l} t^{n+1} \right ) \\
    &=& st \sum_{P \in {\mathrm{R}}} \, \frac{1-s^l}{1-s} \, t^n \, + \,
    2st \sum_{P \in {\mathrm{N}}} \, \frac{1-s^{l}}{1-s} \, t^n \\
    &=& \frac{st}{1-s} \, \left ( R(1,t) - R(s,t) \right ) \, + \, \frac{2st}{1-s} \, \left ( N(1,t) - N(s,t) \right ).
\end{eqnarray*}

\noindent Again, we set $\overline{R}(s,t)=R(1,t)-R(s,t)$; the previous equation can be re-written as:

\begin{equation}\label{nst}
N(s,t) \, \left ( 1+ \frac{2st}{1-s} \right) \, = \, \frac{st}{1-s} \, \overline{R}(s,t) \, + \,
\frac{2st}{1-s} \, N(1,t).
\end{equation}

\noindent As for equation (\ref{rst}), also (\ref{nst}) can be solved using the kernel method. Here the kernel
has a unique solution:

$$ s_1=\frac{1}{1-2t};$$

\noindent replacing $s$ with $s_1$ in $(\ref{nst})$ we have:

$$ \frac{s_1t}{1-s_1} \, \overline{R}(s_1,t) \, + \,
\frac{2s_1t}{1-s_1} \, N(1,t) \, = \, 0, $$

\noindent which leads to:

\begin{equation}\label{n1t}
N(1,t)={\frac {1-7\,t+14\,
{t}^{2}-4\,{t}^{3}}{ \left( 1-2\,t \right)  \left( 1 -4\,t \right) ^{2
}} -{\frac {1-3\,t}{ \left( 1-4\,t \right) ^{3/2}}}}.
\end{equation}

\noindent Replacing the expression of $N(1,t)$ in (\ref{nst}) we can obtain also $N(s,t)$.
Finally we have the generating function:

\begin{equation}\label{fst}
\small
F(s,t)={\frac { \left( 12\,{t}^{3}{s}^{2}-6\,{t}^{2}{s}^{2}+t{s}^{2}-8\,{t}^{
2}s+5\,ts-s+4\,{t}^{2}+1-4\,t \right) ts}{ \left( 1-2\,ts \right)
 \left( 1-s+t{s}^{2} \right)  \left( 1-4\,t \right) ^{2}}}+{\frac {{t
}^{2}s \left( s-2 \right) }{ \left( 1-4\,t \right) ^{3/2} \left( t{s
}^{2}-s+1 \right) }},
\end{equation}

\noindent and the generating function of convex permutominoes according to the size, which gives, after some
simplifications:

\begin{equation}\label{f}
F(1,t)={\frac {2 \, t \left( 1-3\,t \right) }{ \left( 1-4\,t \right) ^{2}}}-{
\frac {t}{ \left( 1-4\,t \right) ^{3/2}}},
\end{equation}

\noindent Starting from (\ref{f}) and performing standard calculations we have the following closed
form for the number $f_n$ of convex permutominoes of size $n$:

\begin{equation}\label{fn}
    f_n= 2 \, (n+3) \, 4^{n-2} \, - \, \frac{n}{2} \, {{2n} \choose {n}} \qquad n\geq 1.
\end{equation}

\noindent The first terms of the sequence are

$$ 1,4,18,84,394,1836,8468, \ldots $$

We remark that while both the left and the right summands of (\ref{fn}) are in \cite{sloane} (sequence A079028
and A002457, respectively), the sequence $\{ f_n \}_{n\geq 0}$ is not present in the Sloane database.

\medskip

Finally we observe that also the number of {\em stack} and {\em directed convex} permutominoes can be easily
obtained from the previous computation.

In fact a stack permutomino can be uniquely represented by a permutomino in class $B$ having the same size. Hence
the generating function of stack permutominoes is given by $B(1,t)$, and then the number of stack permutominoes
of size $n$, as already stated in \cite{fanti}, is equal to $2^n$.

Similarly, a directed convex permutomino can be uniquely represented by a permutomino in class $B$ or one in
class $R$ satisfying {\bf U1}, having the same size. Hence the generating function of directed convex
permutominoes is given by

$$
B(1,t)\, + \, \frac{1}{2} \, R(1,t) \, = \, \frac{t}{1-2t} \, + \, \frac{1}{2} \, \frac {1}{\sqrt {1-4\,t}} - {
\frac{1}{1-2\,t}} \, = \, \frac{1-\sqrt{1-4t}}{2 \, \sqrt{1-4t}},
$$

and then the number of directed convex permutominoes of size $n$, as already stated in \cite{fanti}, is equal to

\begin{equation}\label{dirl}
\frac{1}{2} \, {{2n} \choose {n}}.
\end{equation}

\section{Further work}

In this paper we solve the problem of determining a closed formula for the number of convex permutominoes
with a fixed size. We reach this goal by defining a recursive construction of all the permutominoes
of size $n+1$ starting from those of size $n$, for any $n\geq 1$.

Several problems on the class of permutominoes however remain still open. Below we propose a small
list of the problems which we are interested in, and we would like to tackle in some future work:

\begin{description}
    \item[1.] to give a combinatorial interpretation of the closed formula (\ref{fn}) for the number of
    convex permutominoes. We remark that the expression of $f_n$ resembles the expression (\ref{eq1}) for the
    number of convex polyominoes of fixed semi-perimeter;
    \item[2.] we believe that the ECO construction of convex permutominoes can be extended to the class of
    column-convex permutominoes. One of the problems with this class is that here Proposition \ref{spigoli} does not hold as
    shown in Figure~\ref{colcon}.

    \begin{figure}[htb]
    \begin{center}
    \epsfig{file=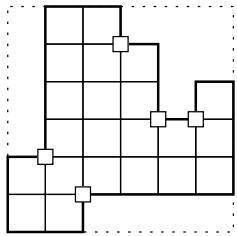} \caption{A column-convex permutominoes: its reentrant points do not satisfy the statement
    of Proposition \ref{spigoli}. \label{colcon}}
    \end{center}
    \end{figure}

    \item[3.] enumerate convex permutominoes according to the {\em area}, i.e. the number of cells of the permutomino.
    The ECO construction we have determined easily leads to a functional equation satisfied by the generating function
    of convex permutominoes according to the area and the size of
    the permutomino; however, then we have not been able to solve this equation.
\end{description}

\end{document}